\documentclass[11pt]{article}

\usepackage{amsmath, amsthm, amsfonts, amssymb}

\newtheorem{thm}{Theorem}[section]
\newtheorem{cor}[thm]{Corollary}
\newtheorem{lem}[thm]{Lemma}
\newtheorem{prop}[thm]{Proposition}
\newtheorem{example}[thm]{Example}
\newtheorem{remarks}[thm]{Remarks}
\newtheorem{defn}[thm]{Definition}

\numberwithin{equation}{section}

\title{\bf Positive curvature property for some hypoelliptic heat
kernels}
\author{Bin Qian  \thanks{Department of Mathematics, Changshu Institute of Technology,
Changshu, Jiangsu 215500,  China, and Institut de Math\'ematiques de
Toulouse, Universit\'e de Toulouse, CNRS 5219.
 E-mail: binqiancn@yahoo.co.cn, binqiancn@gmail.com. Partially supported by the China Scholarship
 council (2007U13020)}
 }\date{}

\newcommand{\rr}{\mathbb{R}}

\def\DD{\mathcal D}

\def\LL{\mathcal L}

\def\vep{\varepsilon}
\def\<{\langle}
\def\>{\rangle}

\def\d"{^{\prime\prime}}

\def\bequ{\begin{equation}}
\def\nequ{\end{equation}}

\def\bdef{\begin{defn}}
\def\ndef{\end{defn}}

\def\bthm{\begin{thm}}
\def\nthm{\end{thm}}

\def\bprop{\begin{prop}}
\def\nprop{\end{prop}}

\def\brmk{\begin{remarks}}
\def\nrmk{\end{remarks}}

\def\bexam{\begin{example}}
\def\nexam{\end{example}}

\def\blem{\begin{lem}}
\def\nlem{\end{lem}}

\def\bcor{\begin{cor}}
\def\ncor{\end{cor}}

\def\bprf{\begin{proof}}
\def\nprf{\end{proof}}

\def\bdes{\begin{description}}
\def\ndes{\end{description}}

\def\part{\partial}
\def\lam{\lambda}
\def\bmax{\begin{matrix}}
\def\nmax{\end{matrix}}

\begin{document}
\maketitle
\begin{abstract}
In this note, we look at some hypoelliptic operators arising from
nilpotent rank 2 Lie algebras. In particular, we concentrate on the
diffusion generated by three Brownian motions and their three L\'evy
areas, which is the simplest extension of the  Laplacian on the
Heisenberg group $\mathbb{H}$.  In order to study contraction
properties of the heat kernel, we show that, as in the case of the
Heisenberg group, the restriction of the sub-Laplace operator acting
on radial functions (which are defined in some precise way in the
core of the paper) satisfies a  non negative Ricci curvature
condition (more precisely a
 $CD(0, \infty)$ inequality), whereas the operator itself does not
 satisfy any $CD(r,\infty)$ inequality. From this we may deduce some
 useful, sharp gradient bounds for the associated heat kernel.

\end{abstract}


\vskip10pt \textbf{2000 MR Subject Classification:} 58J35 43A80

\section{Introduction}
In the study of the long (or small) time behavior ( e.g. gradient
estimates, ergodicity etc.) of simple linear parabolic evolution
equations, one often uses  lower bounds on the Ricci curvature
associated to the generator of the heat kernel, see for example
\cite{Ba97,Ledoux,VSC} and the references therein. But this method
fails in general  in hypoelliptic evolution equations,  since the
Ricci ($\Gamma_2$-) curvature in even the simplest example of the
Heisenberg group can not be bounded below as explained e.g. in
\cite{Juillet,BBBC}. Nevertheless, in the Heisenberg group case,
many properties of the  elliptic case remain true, and we shall
details later some of the most interesting ones.

Let us recall  first some basic facts.

\subsection*{ The elliptic case}
 Let $M$ be a complete Riemannian manifold of dimension
$n$ and let $\LL:=\Delta+\nabla h$, where $\Delta$ is the
Laplace-Beltrami operator. For $t\ge 0$, denote by $P_t$ the heat
semigroup generated by $\LL$ (that is formally $P_{t}= \exp(t\LL)$).
For smooth enough function $f,g$,  one defines (see \cite{Ba97})
$$\aligned
\Gamma(f,g)&=|\nabla f|^2=\frac12(\LL fg-f\LL g-g\LL f),\\
 \Gamma_2(f,f)&=\frac12\big(\LL\Gamma(f,f)-2\Gamma(f,\LL
f)\big)=|\nabla\nabla f|^2+(Ric-\nabla\nabla h)(\nabla f,\nabla
f).\endaligned$$ We have the following well-known proposition, see
Proposition 3.3 in \cite{Ba97}.

\noindent{\bf Proposition A.} For every real $\rho\in \rr$, the
following are equivalent \bdes
\item{(i).} $CD(\rho,\infty)$ holds. That is
$\Gamma_2(f,f)\ge \rho\Gamma(f,f)$.

\item{(ii).} For $t\ge0$, $\Gamma(P_tf,P_tf)\le e^{-2\rho
t}P_t(\Gamma(f,f)).$

\item{(iii).} For $t\ge 0$, $\Gamma(P_tf,P_tf)^{\frac12}\le e^{-\rho
t}P_t(\Gamma(f,f)^{\frac12}).$

\ndes


Moreover, in \cite{Eng06}, Engoulatov obtained the following
gradient estimates for the heat kernels in Riemannian manifolds.

 \noindent {\bf Theorem B. } Let $M$ be a complete Riemannian of
dimension $n$ with Ricci curvature bounded from below, $Ric(M)\ge
-\rho$, $\rho\ge0$. \bdes
\item{(i).} Suppose a non-collapsing condition is satisfies on $M$,
namely, there exist $t_0>0$, and $\nu_0>0$, such that for any $x\in
M$, the volume of the geodesic ball of radius $t_0$ centered at $x$
is not too small, $ Vol(B_x(t_0))\ge \nu_0.$ Then there exist two
constants $C(\rho,n,\nu_0,t_0)$ and $\bar{C}(t_0)>0$, such that
$$
|\nabla \log H(t,x,y)|\le
C(\rho,n,\nu_0,t_0)\left(\frac{d(x,y)}{t}+\frac1{\sqrt{t}}\right),
$$uniformly on $(0, \bar{C}(t_0)]\times M\times M$, where $d(x,y)$
is  the Riemannian distance between $x$ and $y$.

\item{(ii).} Suppose that $M$ has a diameter bounded by $D$, Then
there exists a constant $C(\rho, n)$ such that
$$
|\nabla\log H(t,x,y)|\le
C(\rho,n)\left(\frac{D}{t}+\frac1{\sqrt{t}}+\rho\sqrt{t}\right),
$$
uniformly on $(0,\infty)\times M\times M$.

 \ndes

\subsection*{ The three-dimensional model groups }
 In recent year, some focus has been set on some degenerate
 (hypoelliptic)
situations, where the methods used for the elliptic case do not
apply. Among the simplest examples of such situation are the
three-dimensional groups $\mathfrak{G}$ with Lie algebra
$\mathfrak{g}$, where there is a basis $\{X,Y,Z\}$ of $\mathfrak{g}$
such that \bequ\label{heisenberg-comm} [X,Y]=Z,\ [Z,Y]=\alpha Y,\
[Y,Z]=\alpha X, \nequ where $\alpha\in \rr$.  The analysis reduces
mainly to the thre cases $\alpha=0, \alpha=1, \alpha=-1$.

\bexam[Heisenberg
group, $\alpha=0$] The Heisenberg group can be seen the Euclidean
space $\rr^3$ with a group structure $\circ$,  which is  defined,
 for $\vec{x}=(x,y,z),\vec{y}=(x',y',z')\in\rr^3$, by
$$\vec{x}\circ\vec{y}=\left(x+x',y+y',z+z'+\frac12(xy'-x'y)\right).$$
The left invariant vector fields which are given by
$$\aligned
X(f)&=\lim_{\vep\to
0}\frac{f(\vec{x}\circ(\vep,0,0))-f(\vec{x})}{\vep}
=\left(\part_x-\frac{y}{2}\part_z\right)f, \\
 Y(f)&=\lim_{\vep\to
0}\frac{f(\vec{x}\circ(0,\vep,0))-f(\vec{x})}{\vep}
=\left(\part_y+\frac{x}{2}\part_z\right)f,\\
Z(f)&= \lim_{\vep\to
0}\frac{f(\vec{x}\circ(0,0,\vep))-f(\vec{x})}{\vep} =\part_z
f.\endaligned$$

\noindent The right invariant ones are:
$$\aligned
\hat{X}(f)&=\lim_{\vep\to
0}\frac{f((\vep,0,0)\circ\vec{x})-f(\vec{x})}{\vep}
=\left(\part_x+\frac{y}{2}\part_z\right)f,\\
 \hat{Y}(f)&=\lim_{\vep\to
0}\frac{f((0,\vep,0)\circ\vec{x})-f(\vec{x})}{\vep}
=\left(\part_y-\frac{x}{2}\part_z\right)f,
\endaligned$$

The Lie algebra structure is described by the identities
$[X,Y]=Z,[X,Z]=[Y,Z]=0$. In fact, all group structures satisfying
(\ref{heisenberg-comm}) with $\alpha=0$ can be transformed to  the
case $(\rr^3, \circ)$ by the exponential maps, the vectors fields
$\{X,Y,Z\}$  corresponding to the left ones,  see Lemma 4.1 in
\cite{Gaveau}, see also \cite{BGG}. The natural sublaplacian operator
for this model is $\LL=X^2+Y^2$. In this case, symmetries play an
essential role : they are described by the Lie algebra of the vector
fields that commute with $\LL$.  A basis of this Lie algebra is
$(\hat{X},\hat{Y}, Z)$ and
$\theta=x\part_y-y\part_x$. The last one reflects the rotational
invariance of $\LL$, see \cite{BBBC}.   For this sublaplacian $\LL$,
we have
$$ \Gamma(f,f)=(Xf)^2+(Yf)^2,$$
 and $$
\Gamma_2(f,f)=(X^2f)^2+(Y^2f)^2+\frac12(XYf+YXf)^2+\frac12
(Zf)^2+2\big(XZfYf-YZfXf\big).
$$
The appearance of the mixed term $XZfYf-YZfXf$ prevents the
existence of any constant $\rho\in \rr$ such that $\Gamma_2\ge
\rho\Gamma$. Therefore the methods
used in the elliptic case to prove gradient bounds cannot be used here.
Nevertheless, B. Driver and T. Melcher proved in \cite{DM}, the
existence of a finite positive constant $C_2$ such that
\bequ\label{DM} \forall f\in C^\infty(\mathbb{H},\rr), \ \forall
t\ge 0,\ \Gamma(P_tf,P_tf)\le C_2P_t\Gamma(f,f), \nequ where $P_t$
denotes the associated heat semigroup generated by $\LL$,
$C^\infty(\mathbb{H},\rr)$ is the class of smooth function form
$\mathbb{H}$ to $\rr$ with all partial derivatives of polynomial
growth.  More recently, H. Q. Li \cite{HQLi}
 showed that there exists positive constant $C_1$ such that
\bequ\label{heisenberg-Li} \forall f\in
\mathcal{P}^\infty(\mathbb{H}), \ \forall t\ge 0,\
\Gamma(P_tf,P_tf)^{\frac12}\le
C_1P_t\big(\Gamma(f,f)^{\frac12}\big). \nequ  (See also  D. Bakry et
al. \cite{BBBC} for alternate proofs.) The gradient estimate
(\ref{heisenberg-Li}) is much stronger than  (\ref{DM}), and has
many consequence in terms of functional inequalities for the heat
kernel $P_t$, including Poincar\'e inequalities, Gross logarithmic
Sobolev inequalities, Cheeger type inequalities, and Bobkov type
inequalities, see section 6 in \cite{BBBC}.

Let $p_t$ be the heat
kernel of $P_t$ at $0$ with respect to Lebesgue measures on $\rr^3$.
 In \cite{HQLi}, H. Q. Li
has also pointed out that for $t\ge0, g\in \mathbb{H}$, there exists
 a positive constant $C$ such that \bequ\label{heisenberg-gradient} |\nabla \log
p_t|(g)\le \frac{Cd(g)}{t}, \nequ where $d(g)$ denotes the
Carnot-Carth\'eodory distance (see (\ref{ccdistance})) between $0$
and $g$. This  gradient estimate is sharp and  plays an
important role in the proof of (\ref{heisenberg-Li}).
  \nexam
In the case  $\alpha=1$, the Lie algebra is the one of the  $SU(2)$
Lie group, and
this case has been studied by F. Baudoin and M. Bonnefont in
\cite{Baudoin-Bonnefont}. They show that a modified form of
(\ref{heisenberg-Li}) and (\ref{heisenberg-gradient}) hold. Other
generalizations of Heisenberg group are the so-called   Heisenberg  type group.
They have been studied by H. Q. Li in \cite{HLi1,HLi2},  where he
shows that
(\ref{heisenberg-Li}) and (\ref{heisenberg-gradient}) hold in this
setting. 
In this note, we  shall focus on  a group  that we may call, the
three Brownian motions model.  It can be seen an another typical
simpe example of hypoelliptic operator, but  the structure is more
complex than the Heisenberg (type) groups and the method of H.Q. Li
fails to study the precise gradient bounds in this context.

For this model, we shall first look at the symmetries, that is
characterize all the vector fields  which commute with the
sublaplacian operator $\LL$, see Proposition \ref{linear}. The
infinitesimal rotations are those vector fields which vanish at $0$
and a radial function is a function which vanishes on infinitesimal
rotations. In this case, although the Ricci curvature is everywhere
$-\infty$, refer to \cite{Juillet,BBBC}, we shall prove that the
$\Gamma_2$ curvature is still positive along the radial directions,
as it is the case for the Heisenberg group, see Proposition
\ref{3BM-positive}. As a consequence, the same form of gradient
estimate (\ref{heisenberg-gradient}) holds
 by combining the method developed by F. Baudoin and M. Bonnefont in
\cite{Baudoin-Bonnefont} with the method in \cite{HLi1}. It is worth
recalling that in \cite{BBBQ}, D. Bakry et al. have obtained the
Li-Yau type gradient estimates for the three dimensional model group
by applying $\Gamma_2$-techniques. In our setting, it is easy to see
that this type of gradient estimate  also holds.

\subsection*{The three Brownian motions model}
 The three Brownian motions model $\mathfrak{N}_{3,2}$, see section 4 in \cite{Gaveau}, can be described as the Euclidean
 space $\rr^6$ with a the  following group structure $\circ
 $, which is defined by for $\vec{x}=(x_1,x_2,x_3,y_1,y_2,y_3),\vec{y}=(x_1',x_2',x_3',y_1'
 ,y_2',y_3')\in \rr^6$,
$$
\aligned
(x_1,x_2,x_3,y_1,y_2,y_3)\circ(x_1',x_2',x_3',y_1',&y_2',y_3')=\big(x_1+x_1',x_2+x_2',x_3+x_3',
y_1+y_1'+\frac12(x_2x_3'-x_3x_2'),\\
&\hskip 24pt y_2+y_2'+\frac12(x_3x_1'-x_1x_3'), y_3+y_3'
+\frac12(x_1x_2'-x_2x_1') \big).
\endaligned
$$
For simplification, we make the convention that the index $i\equiv
j\ \mbox{mod} \ 3$, and here we choose $j=1,2,3$. In what follows,
denote $\mathfrak{N}_{3,2}=(\rr^6,\circ)$ be the three Brownian
motions model.

\noindent The three left invariant vector fields which are given,
for $1\le i\le 3$, by
$$
\aligned
X_if&=\lim_{\vep\to0}\frac{f(\vec{x}\circ(\vep_1,\vep_2,\vep_3,0,0,0))
-f(\vec{x})}{\vep}=\left(\part_i-\frac{x_{i+1}}2\hat{\part}_{i+2}+\frac{x_{i+2}}2\hat{\part}_{i+1}\right)f,\\
Y_if&=\lim_{\vep\to0}\frac{f(\vec{x}\circ(0,0,0,\vep_1,\vep_2,\vep_3))
-f(\vec{x})}{\vep}=\hat{\part}_{i}f,
\endaligned
$$
 where $\vep_i=\vep$ and $\vep_j=0$ for $j\neq i$. Here  we use the notation $\hat{\part}_i=\part_{y_i}$.

\noindent The right invariant vector fields which are gives
$\hat{X}_i$, for $1\le i\le 3$, $\vep_i=\vep$ and $\vep_{j}=0$ for
$j\neq i$,
$$
\hat{X}_if=\lim_{\vep\to0}\frac{f((\vep_1,\vep_2,\vep_3,0,0,0)\circ\vec{x})
-f(\vec{x})}{\vep}=\left(\part_i+\frac{x_{i+1}}2\hat{\part}_{i+2}-\frac{x_{i+2}}2\hat{\part}_{i+1}\right)f.
$$
There are no $\hat{Y_i}$'s since in this setting the left and right
multiplications coincide. The Lie algebra structure is described by
the formulae, for $1\le i,j\le
 3$, \bequ\label{three-comm}[X_i,X_{i+1}]=Y_{i+2},\ [X_i,Y_j]=0.\nequ
Similarly for all group structure satisfying (\ref{three-comm}) can
be transformed to  the case  $(\rr^6, \circ)$ via the exponential
maps, the vectors fields are corresponding to the left ones.

 In what
follows, we are interested in the natural sublaplacian  for this model,
which is defined by
$$ \LL=\sum_{i=1}^3X_i^2.
$$
The reason why we call it the three Brownian motions model is that
$\frac12\LL$ is the infinitesimal generator of the Markov process
$\big(\{B_i\}_{1\le i\le
3},\{\frac12\int_0^tB_idB_{i+1}-B_{i+1}dB_{i}\}_{1\le i\le 3}\big)$,
where $\{B_i\}_{1\le i\le3}$ are three  real standard independent
Brownian motions.

For all $t\ge0$, $P_t:=e^{t\LL}$ denotes the associated heat
semigroup generated by $\LL$, $p_t$ the heat kernel of $P_t$ at $0$
with respect to Lebesgue measures on $\rr^6$. For this operator
$\LL$, we have
$$
\Gamma(f,g)=\sum_{i=1}^3X_ifX_ig$$ and
$$
 \Gamma_2(f,f)=\sum_{i,j=1}^3(X_iX_jf)^2-2\sum_{i=1}^3X_if(X_{i+1}Y_{i+2}f
 -Y_{i+1}X_{i+2}f).
$$
Here again the mixed term $\sum_{i=1}^3X_if(X_{i+1}Y_{i+2}f
 -Y_{i+1}X_{i+2}f)$ prevents the existence of any constant $\rho$
 such that
 the curvature dimensional condition $CD(\rho,\infty)$ holds.
 Nevertheless, we have the following Driver-Melcher inequality, see
 \cite{M08},
$$
\Gamma(P_tf,P_tf)\le CP_t(\Gamma(f,f)),
$$
for some positive constant $C$. The constant $C$ here can be
expressed explicitly following  the method in \cite{BBBC}. Also the
optimal reverse local Poincar\'e inequality holds, see Remark 3.3 in
\cite{BBBC}. That is, for any $t\ge0$ and any $f\in
C_c^\infty(\mathfrak{N}_{3,2})$,
$$
t\Gamma(P_tf,P_tf)\le \frac{3}{2}\left(P_t(f^2)-(P_tf)^2\right).
$$
For the H. Q. Li inequality (\ref{heisenberg-Li}), the methods  deeply rely on the precise estimates on
the heat kernel $p_t$ and its differentials (see \cite{HQLi,BBBC}). Up to the author's
knowledge, these precise estimates are not known in the three
Brownian motions model, neither  the H. Q. Li inequality. Nevertheless,
we shall prove that one of the  key gradient estimates (\ref{heisenberg-gradient})
holds, which would be a first step for the proof of the  H. Q. Li
inequality in this context, see Proposition
\ref{3BM-gradient1}.

 The dilation operator in this model is defined by
$\DD:=\frac12\sum_{i}^3x_i\part_i+\sum_{i=1}^3y_i\hat{\part}_{i}$,
and it satisfies \bequ\label{3BM-dilation}[\LL, \DD]=\LL.\nequ For
$t\ge0$, let $T_t=e^{t\DD}$ be the semigroup generated by $\DD$,
that is
$$T_tf(x_1,x_2,x_3,y_1,y_2,y_3)=f\left(\exp{(\frac{t}2)}x_1,\exp{(\frac{t}2)}x_2,\exp{(\frac{t}2)}x_3,\exp{(t)}y_1,
\exp{(t)}y_2,\exp{(t)}y_3\right).$$ From the commutaton relation
(\ref{3BM-dilation}), one deduces, for $t,s\ge0$, $$
P_tT_s=T_sP_{e^st}.$$ Since $0$ is a fixed point of the dilation
group $T_t$, it follows  \bequ\label{scaling} P_t(f)(0)=P_1(T_{\log
t}f)(0).
  \nequ
So it is enough to describe the heat kernel at any tme and any point
to know the operator $P_1(f)(0)$.

The natural distance, induced by the sublaplacian operator $\LL$, is
the Carnot-Carath\'eodory distance $d$.  As usual, it can be defined
from the operator $\LL$ only  by
\bequ\label{ccdistance}
d(g_1,g_2):=\sup_{\{f:\Gamma(f)\le1\}}f(g_1)-f(g_2). \nequ For this
distance, we have the invariant and scaling properties, see
\cite{Gaveau,VSC}.
$$d(g_1,g_2)=d(g_2^{-1}\circ g_1, 0):=d(g_2^{-1}\circ g_1),\
\mbox{and}\  d(\gamma \vec{x},\gamma^2\vec{y})=\gamma d(x,y),$$ for
all $g_1,g_2\in\mathfrak{N}_{3,2}$, $\gamma\in\rr^+$ and
$x=(x_1,x_2,x_3), y=(y_1,y_2,y_3)\in\rr^3$.

\section{Rotation vectors and Radial functions}
In this section, we shall characterize all the vector fields which commute
with  $\LL$. Obviously the right invariant vector
fields $\{\hat{X_i}, Y_i\}_{i=1,2,3}$ commute with $\LL$ since they
commute with $\{X_i\}_{i=1,2,3}$. Like the rotation vector  field
$\theta=x\part_y-y\part_x$ in the Heisenberg group, which commutes
with $\LL$,  there are three rotation vector fields in this case
$$\theta_i=x_{i+1}\part_{i+2}-x_{i+2}\part_{i+1}
+y_{i+1}\hat{\part}_{i+2}-y_{i+2}\hat{\part}_{i+1},\ i=1,2,3. $$ It
is easy to see that $\{\theta_i\}_{i=1,2,3}$ commute with $\LL$ and
we have $[\theta_i,\theta_{i+1}]=\theta_{i+2}$, for $1\le i\le 3$. We
first have the

\bprop\label{linear} The vector fields which  commute with
$\LL$ are the linear combination of the following nine vector fields: the
three right invariant vectors, the three rotations
$\{\theta_i\}_{1\le i\le 3}$, and $\hat{\part}_1,\ \hat{\part}_2,\
\hat{\part}_3$, that is
\bequ\label{prop-linear}\mathcal{T}=\left\{X: X\in
\mbox{span}\{X_i,Y_i,\ 1\le i\le3\},\ [\LL,X]=0\right\}
=\mbox{Linear}\{\hat{X}_i, \theta_i, Y_i, 1\le i\le 3 \}.\nequ

Here
"span" means the  we consider linear combinations of the vector
fields with smoth functions as coefficients, , while "Linear"
 means the coefficients are constants.

\nprop

\bprf We only need to show that the left hand side space in
(\ref{prop-linear}) is contained in the right hand side one. To this
end, for any vector field $X=\sum_{i=1}^3a_iX_i+b_iY_i$ for some
smooth function $a_i,b_i$, satisfies $[\LL,\ X]=0$. For $1\le i\le
3$, denote $Z_i=X_iX_{i+1}+X_{i+1}X_i$, it yields
$X_iX_{i+1}=\frac{Z_{i+2}+Y_{i+2}}{2}$ and
$X_iX_{i+2}=\frac{Z_{i+1}-Y_{i+1}}{2}$. Notice that
$$\aligned\hskip -1pt
[\LL, X]&=\sum_{i=1}^3\Big(\LL a_iX_i
+(\LL b_i+X_{i+1}a_{i+2}-X_{i+2}a_{i+1})Y_i+(X_{i+1}a_{i+2}+X_{i+2}a_{i+1})Z_i\\
&\hskip
12pt+2X_ia_iX_i^2+2X_ib_iX_iY_i+2(X_{i}b_{i+1}-a_{i+2})X_iY_{i+1}+2(X_ib_{i+2}+a_{i+1})X_iY_{i+2}
\Big),\endaligned
$$
thus we have, for $1\le i\le 3$, \begin{gather}
 X_ia_i=X_ib_i=0,\label{comm}\\
\hskip 48pt X_{i+1}b_{i+2}=-X_{i+2}b_{i+1}=a_i,\tag{\ref{comm}$'$}\\
X_{i}a_{i+1}=-X_{i+1}a_{i}.\tag{\ref{comm}$''$} \end{gather} Let us
first  prove the following two claims. \bdes
\item{\bf Claim I:} For $1\le i\le 3$, $a_i$ is independent on $\{y_i,1\le i\le 3\}$ and
linear in $\{x_i,1\le i\le3\}$. To proof the desired result, for
$1\le i\le3$, we have the following commutative property:
$[X_i,Y_{i+1}]=0$, together with (\ref{comm}), it yields
$X_iY_{i+1}b_i=0 $. Since $Y_{i+1}=[X_{i+2}, X_i]$, we can get
$X_i^2X_{i+2}b_i=0$, thus $X_i^2a_{i+1}=0$ by the relation
$X_{i+2}b_i=a_{i+1}$. Similarly we have $X_i^2a_{i+2}=0$. In fine,
together with $X_ia_i=0$,

 $$
X_i^2a_j=0,\ \ i,j=1,2,3. $$

Since
$$
\aligned
 \hskip 12pt [X_1,Y_3]b_2&=0\Rightarrow 2X_1X_2a_3=-X_2X_3a_1,\\
[X_3,Y_1]b_2&=0\Rightarrow2X_3X_1a_2=-X_2X_3a_1,\\
[X_1,Y_2]b_3&=0\Rightarrow 2X_1X_2a_3=-X_3X_1a_2,
\endaligned
$$
we have $$X_1X_2a_3=X_2X_3a_1=X_3X_1a_2=0.$$ Together with the fact
$X_ia_j=-X_ja_i$ by (\ref{comm}$''$), we have $X_iX_ja_k=0$ for
$i,j,k$ all different. Thus we can conclude
\bequ\label{diff-2} X_iX_ja_k=0 \ \ \mbox{for} \ 1\le i,j,k\le
3.\nequ

Note that $Y_1a_1=X_2X_3a_1-X_3X_2a_1, Y_2a_1=X_1^2a_3$ and
$Y_3a_1=-X_1^2a_2$, thanks to (\ref{diff-2}), we get $Y_ia_1=0,
i=1,2,3$. That is $a_1$ is independent on $\{y_i,  i=1,2,3\}$.
Similarly $a_2, a_3$ is independent of $\{y_i,  i=1,2,3\}$. Then
from the definition $X_i$, we have $X_ia_j=\part_ia_j$ for $1\le
i,j\le 3$. With (\ref{diff-2}), we can conclude that $\{a_i,
i=1,2,3\}$ is linear in $\{x_i, i=1,2,3\}$.

  Note that we can also write in the form
  $X=\sum_{i=1}^3a_i\part_i+c_i\hat{\part}_i$,
where $c_i=b_i+\frac12(a_{i+2}x_{i+1}-a_{i+1}x_{i+2}), \ 1\le i\le
3.$ Then we can conclude

\item{\bf Claim II:} $\{c_i,i=1,2,3\}$ is linear in $\{x_i, y_i,\  1\le i\le 3\}$.
By Claim I, $a_i$ is independent of $y_j$, together with the fact
$X_ia_i=0$, we have the equation (\ref{comm}$'$)
 is equivalent to \bequ\label{diff-11}\aligned
\frac12a_i&=X_{i+1}c_{i+2}+\frac12x_{i+1}\part_{i+1}a_i\\
&=-X_{i+2}c_{i+1}+\frac12x_{i+2}\part_{i+2}a_i.
\endaligned\nequ
And $X_ib_i=0$ is equivalent to \bequ\label{diff-12}
X_ic_i=\frac12(x_{i+1}\part_ia_{i+2}-x_{i+2}\part_ia_{i+1}). \nequ
Using (\ref{comm})-(\ref{diff-12}), the relations
$[X_i,X_{i+1}]=Y_{i+2}$ and Claim I, through computation, we have
\bequ\label{diff-13}Y_ic_j=\part_ia_j,\ \mbox{for} \ i,j=1,2,3.\nequ

Since $a_i$ is linear in $x_i$, we can conclude that $c_j$ has no
second order terms in $\{x_i,y_i,1\le i\le3\}$.  By the definition
of $X_i$ and (\ref{diff-12}) and (\ref{diff-13}), we have
$$
\aligned
\part_ic_i&=X_ic_i+\frac{x_{i+1}}2Y_{i+2}c_i-\frac{x_{i+2}}2Y_{i+1}c_i\\
&=\frac12(x_{i+1}\part_ia_{i+2}-x_{i+2}\part_ia_{i+1})-\frac{x_{i+1}}2\part_ia_{i+2}+\frac{x_{i+2}}2\part_ia_{i+1}\\
&=0.\endaligned$$ By (\ref{diff-11}) and (\ref{diff-13}),
\bequ\label{diff-14} \aligned
\part_ic_{i+1}&=X_ic_{i+1}+\frac{x_{i+1}}2Y_{i+2}c_{i+1}-\frac{x_{i+2}}2Y_{i+1}c_{i+1}\\
&=\frac12a_{i+2}-\frac12x_i\part_ia_{i+2}-\frac12x_{i+1}\part_{i+1}a_{i+2},\endaligned
\nequ similarly, \bequ\label{diff-15}
\part_ic_{i+2}=-\frac12a_{i+1}+\frac12x_i\part_ia_{i+1}+\frac12x_{i+2}\part_{i+2}a_{i+1}.
\nequ By Claim I, (\ref{diff-13})-(\ref{diff-15}) and
$\part_ia_i=0$, we can conclude that for $1\le i,j\le 3$,
$\part_ic_j$ is constant. Thus we complete to proof Claim II. \ndes
 By the above two claims and
$\part_ia_i=0$, we
can assume $$ \aligned a_i&=A_{i,i+1}x_{i+1}+A_{i,{i+2}}x_{i+2}+B_i,\\
\endaligned
$$
where $A_{i,j}=-A_{j,i}$, $B_i$ are constants, then we have, by
(\ref{diff-13})-(\ref{diff-15}),

$$
\aligned
c_i&=\frac12(B_{i+1}x_{i+2}-B_{i+2}x_{i+1})+A_{i,i+1}y_{i+1}
-A_{{i+2},i}y_{i+2}+D_i,\endaligned
$$
where $D_i$ are constants.

If we choose $B_i=1$ (or respectively $D_i=1$, $A_{i,i+2}=1$) and
the other constants  $0$, we get $X=\hat{X}_i$ (or respectively
$Y_i$, $\theta_i$). Thus we complete the proof.

\nprf In the Heisenberg group, the radial functions $f$ can be
characterized by $\theta f=0$. Here in our setting, as an extension
of such characterization, we can give a definition of radial
functions. \bdef A smooth enough function $f$ is called radial if
and only if for
 $1\le i\le 3$, $\theta_if=0$. \ndef

 (Notice that here the vector fields $\theta_{i}$ are the commuting
 vector fields which vanish in $0$.)

 \brmk
Note that the heat kernel $(p_t)_{t\ge0}$ is radial. The reason is
that for any function $f$, $1\le i\le3$, $\theta_i f(0)=0$ and
$\{\theta_i\}_{1\le i\le3}$ commute with $\LL$, whence they commute
with the semigroup $P_t=e^{t\LL}$. Hence, for any function $f$, one
has $P_{t}\theta_{i} f=0$, which, taking the adjoint of $\theta_{i}$
under the Lebesgue measure, which is $-\theta_{i}$, shows that for the density $p_{t}$ of the
heat kernel at $0$, one has $\theta_{i} p_{t}=0$. This explains why
any information about the radial functions in turns give information
on the heat kernel itself.
 \nrmk
\brmk\label{radial-para}  For any radial function $f$, there exist
some function $g$ such that $f(\vec{x},\vec{y})=g(r_1,r_2,z)$, where
$r_1=\sum_{i=1}^3x_i^2,\ r_2=\sum_{i=1}^3y_i^2,\
z=\sum_{i=1}^3x_iy_i$. Indeed, by the definition, for $1\le i\le 3$,
$\theta_if=0$, then we have $f=f(U\vec{x},U\vec{y})$, where $U$ is
arbitrary linear orthogonal transformation on $\rr^3$, which
satisfying $U^*U=UU^*=1$. Hence
$f=\bar{f}(|\vec{x}|,|\vec{y}|,\<\vec{x},\vec{y}\>)$, for some
function $\bar{f}$. Here is another way, by the transformation
$\theta_i$, we will directly get
$$
f(\vec{x},\vec{y})=f(\sqrt{r_1},0,0,\frac{z}{\sqrt{r_1}},0,
\sqrt{r_2-z^2/r_1}).
$$
 \nrmk

\section{$\Gamma_2$ curvature}
Recall that we can't find a constant $\rho\in\rr$ such that
$\Gamma_2\ge \rho\Gamma$ because of the appearance of the items
$X_ifX_jY_{k}f$. In other words, the  Ricci curvature is everywhere
$-\infty$. Nevertheless we shall prove $\Gamma_2$ curvature is
positive on the radial functions.

\bprop\label{3BM-positive} For any smooth radial function $f$, we
have
$$ \Gamma_2(f,f)\ge0.
$$ \nprop

Here we will give two different proofs.  The first one is that we
shall use directly the three equations asserting that a function is
radial. Then, applying the vector fields $\{X_j\}_{1\le j\le 3}$ on
these equations, we get nine equations in hand. It follows that we
can get the exact expressions of $\{X_iY_jf\}_{1\le i,j\le 3}$ in
terms of $X_{i}X_{j}f$ and also first order terms. (In fact, we
adapt the mathematical software MAPLE to do it). Then we substitute
them into the formal expression of $\Gamma_2$, and we find that
$\Gamma_2$ can also be expressed in a functional non negative
quadratic form.

The second way is that by the Remark \ref{radial-para}, we have an
expression of the sublaplacian operator acting on radial functions
directly through a good parametrization, say $r_1,r_2,z$. Through
 computation, we can obtain the exact expression of $\Gamma_2$ curvature and find again
that $\Gamma_2$ can be expressed in a functional non negative
quadratic form, thus we are done.

\bprf[The first proof]
 A radial function $f$ satisfies $\theta_if=0$, which is equivalent
 to
 say that
 $$ x_{i+1}X_{i+2}f-x_{i+2}X_{i+1}f
 -\frac{x_{i+1}^2+x_{i+2}^2}{2}Y_if
 +\frac{x_ix_{i+1}-2y_{i+2}}{2}Y_{i+1}f
 +\frac{x_ix_{i+2}+2y_{i+1}}{2}Y_{i+2}f=0.
$$
Differentiating the above equations at the directions $\{X_j\}_{1\le
j\le 3}$, with the commutative  relations $[X_{i},X_{i+1}]=Y_{i+2}$,
we get the nine differential equations, for $1\le i\le 3$,
$$\aligned
&x_{i+1}X_{i+2}X_if-x_{i+2}X_{i+1}X_if-\frac{x_{i+1}^2+x_{i+2}^2}{2}X_iY_if
\\
&\hskip 100pt+\frac{x_{i}x_{i+1}-2y_{i+2}}{2}X_iY_{i+1}f+\frac{x_ix_{i+2}+2y_{i+1}}{2}X_1Y_3f=0,\\
&X_{i+2}f+x_{i+1}X_{i+2}X_{i+1}f-x_{i+2}X_{i+1}^2f
-\frac{x_{i+1}^2+x_{i+2}^2}{2}X_{i+1}Y_if \\
&\hskip 100pt+\frac{x_ix_{i+1}
-2y_{i+2}}{2}X_{i+1}Y_{i+1}f+\frac{x_ix_{i+2}
+2y_{i+1}}{2}X_{i+1}Y_{i+2}f=0,\\
&-X_{i+1}f+x_{i+1}X_{i+2}^2f-x_{i+2}X_{i+1}X_{i+2}f
-\frac{x_{i+1}^2+x_{i+2}^2}{2}X_{i+2}Y_if\\
&\hskip 100pt+\frac{x_ix_{i+1}
-2y_{i+2}}{2}X_{i+2}Y_{i+1}f+\frac{x_ix_{i+2}
+2y_{i+1}}{2}X_{i+2}Y_{i+2}f=0.
\endaligned
$$

For simplificity, we will use the following notations, for $1\le
i\le3$,
$$\alpha_{i+1}:=x_{i+1}X_if-x_iX_{i+1}f,\
\beta_{i+1}:=y_{i+1}X_if-y_iX_{i+1}f,$$
$$\gamma_i:=x_{i}y_{i+1}-x_{i+1}y_i,\
\eta_i:=x_{i}y_i+x_{i+1}y_{i+1},
$$
$$  |x|_{i}^2:=x_i^2+x_{i+1}^2,\
A_i:=\gamma_2X_1X_{i+1}f+\gamma_3X_2X_{i+1}f +\gamma_1X_3X_{i+1}f,$$
and $$|h|^2:=\sum_{i=1}^3|h_i|^2, \ \mbox{ for }
h=(h_1,h_2,h_3)\in\rr^3.$$

 From the above nine differential equations, we can get, for
$1\le i\le 3$,
$$\aligned
X_iY_{i+1}f&=-\frac{1}{2|\gamma|^2}
\Big((x_ix_{i+1}|x|^2+2y_{i+2}|x|_i^2-2x_{i+2}\eta_i+4y_{i}y_{i+1})\cdot(x_{i+2}X_{i+1}X_if-x_{i+1}X_{i+2}X_if)\\
&\hskip 24pt
+(x_{i+1}x_{i+2}|x|^2-2y_i|x|_{i+1}^2+2x_i\eta_{i+1}+4y_{i+1}y_{i+2})
\cdot(x_{i+1}X_i^2f-x_iX_{i+1}X_if-X_{i+1}f)\\
&\hskip 24pt
+(x_{i+1}^2|x|^2+4y_{i+1}^2)\cdot(x_iX_{i+2}X_{i}f-x_{i+2}X_i^2f+X_{i+2}f)\Big),\endaligned
$$
and
$$
\aligned X_iY_{i+2}f&=-\frac{1}{2|\gamma|^2}
\Big((x_ix_{i+2}|x|^2-2y_{i+1}|x|_{i+2}^2+2x_{i+1}\eta_{i+2}+4y_{i+2}y_{i})
\cdot(x_{i+2}X_{i+1}X_if-x_{i+1}X_{i+2}X_if)\\
&\hskip 24pt
+(x_{i+1}x_{i+2}|x|^2+2y_i|x|_{i+1}^2-2x_i\eta_{i+1}+4y_{i+1}y_{i+2}
)\cdot(x_{i}X_{i+2}X_{i}f-x_{i+2}X_i^2f+X_{i+2}f)\\
&\hskip
24pt+(x_{i+2}^2|x|^2+4y_{i+2}^2)\cdot(x_{i+1}X_i^2f-x_iX_{i+1}X_if-X_{i+1}f)
\Big).
 \endaligned
$$
Note that
$$
 \Gamma_2(f,f)=\sum_{i,j=1}^3(X_iX_jf)^2-2\sum_{i=1}^3X_if(X_{i+1}Y_{i+2}f
 -Y_{i+1}X_{i+2}f).
$$

With the exact expressions of $X_iY_jf$ in hand, through
calculation, we have
$$\aligned
F:&=2|\gamma|^2\cdot\sum_{i=1}^3X_if(X_{i+1}Y_{i+2}f-X_{i+2}Y_{i+1}f)\\
&=2\sum_{i,j=1}^3\Big(X_{i+j}f\big(2y_{i+j+1}\gamma_{i+1}
-x_{i+j+1}y_i|x|_{i+1}^2+x_ix_{i+j+1}\eta_{i+1}\big)
\\
&\hskip
48pt-X_{i+j+1}f\big(2y_{i+j}\gamma_{i+1}-x_{i+j}y_i|x|_{i+1}^2
+x_ix_{i+j}\eta_{i+1}
\big)\Big)X_{i}X_{i+j+2}f-|x|^2\cdot|\alpha|^2-4|\beta|^2\\
&=2\sum_{i,j=1}^3\Big(2\beta_{i+j}\gamma_{i+1}+\alpha_{i+j}
\gamma_{i}x_{i+1}
-\alpha_{i+j}\gamma_{i+2}x_{i+2}\Big)X_iX_{i+j+1}f-|x|^2\cdot|\alpha|^2-4|\beta|^2.
\endaligned
$$
Rearrange the items, we have
$$
\aligned
F=2\sum_{i,j=1}^3\alpha_ix_{i+j+1}\big(\gamma_{i+j}X_{i+j}X_{i+1}f
-\gamma_{i+j+1}X_{i+j+2}X_{i+1}f\big)+4\sum_{i=1}^3\beta_iA_i-|x|^2\cdot|\alpha|^2-4|\beta|^2.
\endaligned
$$
Notice that
$$
|\gamma|^2\cdot\sum_{i,j=1}^3(X_iX_jf)^2-\sum_{i=1}^3A_i^2
=\sum_{i,j=1}^3\big(\gamma_{i+j}X_{i+j}X_{i+1}f
-\gamma_{i+j+1}X_{i+j+2}X_{i+1}f\big)^2,$$
 it follows that  $$ \aligned
|\gamma|^2\cdot\Gamma_2&=|\gamma|^2\cdot\sum_{i,j=1}^3(X_iX_jf)^2-F\\
&=\sum_{i=1}^3(2\beta_i-A_i)^2+\sum_{i,j=1}^3\big(\gamma_{i+j}X_{i+j}X_{i+1}f
-\gamma_{i+j+1}X_{i+j+2}X_{i+1}f\big)^2\\
&\hskip 24pt -2\sum_{i,j=1}^3\alpha_ix_{i+j+1}\big(\gamma_{i+j}X_{i+j}X_{i+1}f-\gamma_{i+j+1}X_{i+j+2}X_{i+1}f\big)+|x|^2\cdot|\alpha|^2\\
&=\sum_{i=1}^3(2\beta_i-A_i)^2
+\sum_{i,j=1}^3(\gamma_{i+j}X_{i+j}X_{i+1}f-\gamma_{i+j+1}X_{i+j+2}X_{i+1}f
-\alpha_ix_{i+j+1})^2.
\endaligned
$$
Hence we complete the proof.

\nprf

\bprf[The second proof ]

Denote $r_1=\sum_{i=1}^3x_i^2,\ r_2=\sum_{i=1}^3y_i^2,\
z=\sum_{i=1}^3x_iy_i$. Through calculation, we have
$$
\LL r_1=6,\ \LL r_2=r_1,\ \LL z=0,
$$
and
$$\aligned
\Gamma(r_1,r_1)&=4r_1,\ \Gamma(r_2,r_2)=r_1r_2-z^2,\
\Gamma(z,z)=r_2,\\
\Gamma(r_1,r_2)&=0, \ \Gamma(r_1,z)=2z,\ \Gamma(r_2,z)=0.
\endaligned$$

For any radial functions $f,g$ depend only on $r_1,r_2,z$, we have
$$\aligned
\LL f(r_1,r_2,z)&=\part_{r_1}f\LL r_1+\part_{r_2}f\LL
r_2+\part_zf\LL
z+\part_{r_1r_1}^2f\Gamma(r_1,r_1)+\part_{r_2r_2}^2f\Gamma(r_2,r_2)+\part_{zz}^2f\Gamma(z,z)\\
&\hskip
24pt+2\part_{r_1z}^2f\Gamma(r_1,z)+2\part_{r_2z}^2f\Gamma(r_2,z)+2\part_{r_1r_2}^2f\Gamma(r_1,r_2)\\
&=4r_1\part_{r_1r_1}^2f+(r_1r_2-z^2)\part_{r_2r_2}^2f+r_2\part_{zz}f+4z\part_{r_1z}^2f
+6\part_{r_1}f+r_1\part_{r_2}f\\
&:=\hat{\LL} f,
\endaligned
$$
where $\hat{\LL}$ has the following expression
$$
\hat{\LL}
f=4r_1f_{11}+(r_1r_2-z^2)f_{22}+r_2f_{zz}+4zf_{1z}+6f_1+r_1f_2.
$$
Hence for any radial function $f=f(r_1,r_2,z),g=g(r_1,r_2,z)$,
$$\aligned
\Gamma(f,g)&:=\frac12(\LL(fg)-f\LL g-g\LL f)\\
&=\frac12(\hat{\LL}(fg)-f\hat{\LL} g-g\hat{\LL} f)\\
&=\hat{\Gamma}(f,g),\\
\endaligned$$
and also

$$
\aligned
\Gamma_2(f,f)&:=\frac12(\LL\Gamma(f,f)-2\Gamma(f,\LL f))\\
&=\frac12(\hat{\LL}\hat{\Gamma}(f,f)-2\hat{\Gamma}(f,\hat{\LL} f))\\
&=\hat{\Gamma}_2(f,f).
\endaligned
$$
Through direct calculation, we have
$$\aligned
\Gamma(f,g)&=\hat{\Gamma}(f,g)\\
&=4r_1f_1g_1+(r_1r_2-z^2)f_2g_2+r_2f_zg_z+2zf_1g_z+2zf_zg_1\endaligned
$$
and
$$
\aligned \Gamma_2(f,f)&=\hat{\Gamma}_2(f,f)\\
&=16r_1^2f_{11}^2+16r_1f_1f_{11}+8r_1(r_1r_2-z^2)f_{12}^2
+8(r_1r_2-z^2)f_2f_{12}+8(r_1r_2+z^2)f_{1z}^2\\
&\hskip
24pt+32r_1zf_{11}f_{1z}+r_1(r_1r_2-z^2)f_2f_{22}+(r_1r_2-z^2)^2f_{22}^2+
2(r_1r_2-z^2)f_zf_{2z}
\\
&\hskip 24pt+8z(r_1r_2-z^2)f_{12}f_{2z}+(2r_2+\frac{r_1^2}{2})f_2^2
+2r_2(r_1r_2-z^2)f_{2z}^2+r_2^2f_{zz}^2
+4r_2f_1f_{zz}\\
&\hskip 24pt+8r_2zf_{1z}f_{zz}+16zf_1 f_{1z}+
8z^2f_{11}f_{zz}+12f_{1}^2+\frac12r_1f_z^2-4(r_1r_2-z^2)f_1f_{22}
\\
&\hskip 24pt-4r_1f_1f_2-(r_1r_2-z^2)f_2f_{zz} -2zf_2f_z.\endaligned
$$

By  careful study, we can express the above into a functional
quadratic form.
$$
\aligned\Gamma_2(f,f)
&=\left((r_1r_2-z^2)f_{22}+\frac{r_1}{2}f_2-2f_1\right)^2+8r_1(r_1r_2-z^2)\left(f_{12}
+ \frac{f_2}{2r_1}+\frac{zf_{2z}}{2r_1}\right)^2\\
&\hskip 24pt
+\frac{2}{r_1}\left((r_1r_2-z^2)f_{2z}+\frac{r_1}{2}f_z-zf_2\right)^2
+\left(\frac{r_1}2f_2-2f_1-\frac{r_1r_2-z^2}{r_1}f_{zz}\right)^2
\\&\hskip 24pt +4\left(f_1+\frac{z^2}{2r_1}f_{zz}+2r_1f_{11}+2zf_{1z}\right)^2
+2(r_1r_2-z^2)\left(\frac{z}{r_1}f_{zz}+2f_{1z}\right)^2.\endaligned
$$
Hence the desire result follows. \nprf

\section{Gradient bounds for the heat kernels}
As done in \cite{BBBQ}, we have the following Li-Yau type inequality
holds. \bprop\label{prop-LY} There exist positive constants
$C_1,C_2,C_3$ such that for any positive function $f$, if $u=\log
P_tf$, we have
$$
\part_t u\ge C_1\Gamma(u)+C_2t\sum_{i=1}^3|Y_iu|^2-\frac{C_3}{t}.
$$
 \nprop
\bprf Here we briefly proof it for the readers' convenience. Notice
that for all $\lam>0$,
\begin{align}
\Gamma_2(f,f)&=\sum_{i=1}^3(X_i^2f)^2+\frac12\sum_{i=1}^3(Y_if)^2
+2\sum_{i=1}^3\big(D_{i,i+1}(f)\big)^2+2\sum_{i=1}^3\big(X_ifX_{i+2}Y_{i+1}f-X_{i+2}fX_{i}Y_{i+1}f\big)\nonumber\\
&\ge \frac{1}{3}(\LL
f)^2+\frac12\sum_{i=1}^3(Y_if)^2-4\sqrt{\Gamma(f)}\cdot\sqrt{\sum_{i=1}^3\Gamma(Y_if)}\label{gamma_2-inequ}\\
&\ge \frac{1}{3}(\LL
f)^2+\frac12\sum_{i=1}^3(Y_if)^2-\frac{4}{\lam}\Gamma(f)-\lam\sum_{i=1}^3\Gamma(Y_if),\nonumber
\end{align}
where $D_{i,i+1}=\frac12(X_iX_{i+1}+X_{i+1}X_i)$, and the last two
inequalities follow from the Cauchy-Schwartz inequality.

Set $f_s=P_{t-s}f, u_s=\log f_s$, following \cite{BBBQ}, let $$
\Phi_1(s)=P_s\big(f_s\Gamma(u_s,u_s)\big),\
\Phi_2(s)=P_s\left(f_s\sum_{i=1}^3(Y_iu_s)^2\right),
$$
we have $$ \Phi_1'(s)=2P_s\big(f_s\Gamma_2(u_s,u_s)\big),\
\Phi_2'(s)=2P_s\left(f_s\sum_{i=1}^3\Gamma(Y_iu_s)\right).
 $$
Combining (\ref{gamma_2-inequ}) and $$ (\LL u_s)^2\ge 2\gamma \LL
u_s-\gamma^2,\ \LL u_s=\frac{\LL f_s}{f_s}-\Gamma(u_s),
$$
we have $$ \Phi_1'(s)\ge
\left(-\frac{4}{\lam}-\frac{4\gamma}{3}\right)\Phi_1(s)+\Phi_2(s)-\lam
\Phi_2'(s)+\frac{4\gamma}{3}\LL P_tf- \frac{2\gamma^2}{3}P_tf.
$$
Denote $a,b$ are positive functions defined on $[0,t)$, with $b$ is
decreasing, we have
$$
(a(s)\Phi_1(s)+b(s)\Phi_2(s))'\ge
\left(a'-\frac{4a}{\lam}-\frac{4a\gamma}{3}\right)\Phi_1(s)+(a+b')\Phi_2(s)
+(b-\lam a)\Phi_2'+\frac{4\gamma a}{3}\LL P_tf-
\frac{2\gamma^2a}{3}P_tf.
$$
By choosing $$ a=-b',\ \lam=-\frac{b}{b'},\
\gamma=\frac{3b''}{4b'}+\frac{3b'}{b},
$$and then choose $b(s)=(t-s)^{\alpha}$, for some $\alpha>2$,
 integrating the above differential inequality from $0$ to $t$,
the desired result follows.
 \nprf
As a consequence, we have the following Harnack inequality: There
exist positive constants $A_1$, $A_2$, for $t_2>t_1>0$, and
$g_1,g_2\in \mathfrak{N}_{3,2}$, \bequ\label{3BM-harnack}
\frac{p_{t_1}(g_1)}{p_{t_2}(g_2)}\le
\left(\frac{t_2}{t_1}\right)^{A_1}e^{A_2\frac{d^2(g_1,g_2)}{t_2-t_1}}.
\nequ

Here is an analogue result of Theorem B in the three Brownian
motions model.

\bprop\label{3BM-gradient1} There exists a constant $C>0$ such that
for $t>0$, $g=(x,y)\in\mathfrak{N}_{3,2}$,
$$
\sqrt{\Gamma(\log p_t)(g)}\le
\frac{Cd(g)}{t},
$$
where $p_t(g)$ denotes the density of $P_t$ at $0$ and $d(g)$
denotes the Carnot-Carath\'eodory distance between $0$ and $g$.
\nprop \bprf
Following \cite{Baudoin-Bonnefont}, for $0<s<t$, let
$\Phi(s)=P_s\big(p_{t-s}\log p_{t-s}\big)$, we have
$$ \Phi'(s)=P_s\big(p_{t-s}\Gamma(\log p_{t-s})\big),\
\Phi''(s)=2P_s\big(p_{t-s}\Gamma_2(\log p_{t-s})\big).
$$
By Proposition \ref{3BM-positive}, $\Phi''$ is positive, whence
$\Phi'$ is non-desceasing, thus $$\int_0^{\frac{t}2}\Phi'(s)ds\ge
\frac{t}2\Phi'(0).
$$
That is
$$
p_t\Gamma(\log p_t)\le \frac{2}{t}\big(P_{t/2}(p_{t/2}\log
p_{t/2})-p_t\log p_t\big).
$$
The right hand side can be bounded by applying the above Harnack
inequality (\ref{3BM-harnack}) and the basic fact $p_{t/2}(g)\le
p_{t/2}(0)$, for all $g\in \mathfrak{N}_{3,2}$. We have
$$
\sqrt{\Gamma(\log p_t)(g)}\le
C\left(\frac{d(g)}{t}+\frac{1}{\sqrt{t}}\right).
$$
In particular,
$$
\sqrt{\Gamma(\log p_1)(g)}\le C\left(d(g)+1\right).
$$
If $d(g)\ge1$, it is trivial to get the desired result.

 Note that for $x=(x_1,x_2,x_3),\ y=(y_1,y_2,y_3)$, see P. 125, Theorem
 1 in  \cite{Gaveau},
$$
\aligned p_1(x,y)&=-(2\pi)^{-\frac{15}2}\int_{\rr^3}\exp{\left(-i
y\cdot \alpha\right)}\frac{|\alpha|}2\left(\sinh
\frac{|\alpha|}2\right)^{-1}\\
&\hskip 30
pt\cdot\exp{\frac12\left\{-|x|^2-\frac{xA^2x^t}{|\alpha|^2}\left(1-\frac{|\alpha|
}2\coth \frac{|\alpha|}2\right)\right\}} \prod_{k=1}^3d\alpha_{k}\\
&\stackrel{(*)}{=}-(2\pi)^{-\frac{15}2}\int_{\rr^3}\exp{\left(-i
y\cdot \alpha\right)}\frac{|\alpha|}2\left(\sinh
\frac{|\alpha|}2\right)^{-1}\\
&\hskip 30
pt\cdot\exp{\frac12\left\{-|x|^2+(|x|^2-(x\cdot\tilde{\alpha})^2)\left(1-\frac{|\alpha|
}2\coth \frac{|\alpha|}2\right)\right\}} \prod_{k=1}^3d\alpha_{k},
 \endaligned$$
where $\alpha=(\alpha_1,\alpha_2,\alpha_3),
|\alpha|^2=\sum_{k=1}^3\alpha_k^2,\
\tilde{\alpha}=\frac{1}{|\alpha|}(\alpha_1,\alpha_2,\alpha_3) $ and
$
A=\left(\begin{smallmatrix}0&\alpha_1&-\alpha_2\\-\alpha_1&0&\alpha_3\\\alpha_2&-\alpha_{3}&0
\end{smallmatrix}\right)$. $(*)$ follows from
$xA^2x^t=-|\alpha|^2(|x|^2-(x\cdot\alpha)^2)$. Notice
$$\aligned
\part_jp_1(x,y)&=(2\pi)^{-\frac{15}2}\int_{\rr^3}\exp{\left(-i
y\cdot \alpha\right)}|\alpha|\left(\sinh
\frac{|\alpha|}2\right)^{-1}
\cdot\left(-x_j+(x_j-x\cdot\tilde{\alpha}\tilde{\alpha}_j)(1-\frac{|\alpha|}2\coth\frac{|\alpha|}{2})\right)\\
&\hskip 30
pt\cdot\exp{\frac12\left\{-|x|^2+(|x|^2-(x\cdot\tilde{\alpha})^2)\left(1-\frac{|\alpha|
}2\coth \frac{|\alpha|}2\right)\right\}} \prod_{k=1}^3d\alpha_{k}
\endaligned$$
and

$$\aligned
\hat{\part}_ip_1(x,y)&=(2\pi)^{-\frac{15}2}i\int_{\rr^3}\exp{\left(-i
y\cdot \alpha\right)}\frac{|\alpha|}2\left(\sinh
\frac{|\alpha|}2\right)^{-1}\cdot \alpha_i \\
&\hskip 30
pt\cdot\exp{\frac12\left\{-|x|^2+(|x|^2-(x\cdot\tilde{\alpha})^2)\left(1-\frac{|\alpha|
}2\coth \frac{|\alpha|}2\right)\right\}} \prod_{k=1}^3d\alpha_{k}.
\endaligned$$
Let
$$W_1=\int_{\rr^3}|\alpha|\left(\sinh\frac{|\alpha|}{2}\right)^{-1}
\prod_{k=1}^3d\alpha_k,\
W_2=\int_{\rr^3}\left(\sinh\frac{|\alpha|}{2}\right)^{-1}\cdot
|\alpha|^2\coth\frac{|\alpha|}{2}\prod_{k=1}^3d\alpha_k,$$ obviously
we have $W_1,W_2$ are bounded. For $g=(x,y)\in \mathfrak{N}_{3,2}$,
satisfying $d(g)\le 1$, by the basic fact that $|x|\le d(g)\le 1$
(see \cite{NSW,VSC}, in fact we can easily proof it on the nilpotent
groups.), we have
$$
\sqrt{\Gamma(p)(g)}\le C_1|x|(W_1+W_2)\le C_2|x|\le Cd(g).
$$
The desired result follows by the time scaling property
(\ref{scaling}).

\nprf

Notice that for any  radial function $f$, $P_tf$ is also radial
since all $\theta_i$ commute with  $P_t$. Thanks to Proposition
\ref{3BM-positive}, we have H. Q. Li inequality, LSI inequality,
isoperimetric inequalities etc. hold for the semigroup restricted on
the radial functions, see \cite{BBBC,Ba97}.  We state them in the
following proposition. \bprop For any compactly supported smooth,
radial function $f$, for any $t\ge0$, $g\in\mathfrak{N}_{3,2}$,
\bdes
\item{(i) H. Q. Li inequality.}  $\Gamma(P_tf,P_tf)^{\frac12}(g)\le
P_t(\Gamma(f,f)^{\frac12})(g)$.
\item{(ii) LSI inequality.}  $P_t(f\log f)(g)-P_t(f)\log P_t(f)(g)\le
tP_t\left(\frac{\Gamma(f,f)}{f}\right)(g).$
\item{(iii) Isoperimetric inequality.}  $ P_t(|f-P_t(f)(g)|)(g)\le 4\sqrt{t}P_t(\Gamma(f)^{\frac12})(g).$
\ndes \nprop

\vskip 20pt {\bf Discussion:} Here we have shown that H. Q. Li
inequality holds for the radial functions. For the general
functions, it is still open. Following the viewpoint of \cite{HQLi},
also \cite{BBBC}, one key point to proof H. Q. Li inequality is  the
precise lower and upper bounds for the associated heat kernel. But
in our setting, this estimates is unknown, at least the methods in
\cite{HQLi}-\cite{HLi2} are not applicable.  This precise estimates
are also essential to proof the cheeger type inequality, see Lemma
5.1 in \cite{BBBC} and proof the constant coefficient is bounded in
the complex quasi-communication method, see Proposition 5.5 in
\cite{BBBC}.


\end{document}